%% file: newproof.tex
\documentclass[nameyear,final]{article}
\usepackage[english]{babel} 
\usepackage{amsmath,amssymb,amsfonts} 
\usepackage{theorem}
\usepackage{color}
\usepackage[final]{graphicx} 
\usepackage[longnamesfirst]{natbib} 
\usepackage{xargs} 
\usepackage{showkeys} 
\usepackage{enumerate} 

\newtheorem{theo}{Theorem}
\newtheorem{coro}{Corollary}
\newtheorem{lemma}{Lemma}
\newtheorem{propo}{Proposition}
\newcommand{\dd}{\mathrm{d}}
\newcommand{\eps}{\varepsilon}

\newcommand{\cvg}[2][{}]{\overset{#1}{\underset{#2}{\longrightarrow}}}
\newcommand{\eqdistr}{\overset{\mathcal{D}}{=}}
\newcommand{\Stable}[1][t]{\mathcal{S}_{#1}}
\newcommand{\Peaks}[1][t]{\mathcal{H}_{#1}}
\newcommand{\Well}[1][t]{\mathcal{W}_{#1}}
\newcommandx*{\Typical}[2][1,2={t,\eps}]{\Gamma^{#1}_{#2}}
\newcommand{\Neigh}[1][\eps\log t]{D_{#1}}
\newcommand{\SET}[1]{\left\{{#1}\right\}}
\newcommand{\NN}{\mathbb{N}}
\newcommand{\ZZ}{\mathbb{Z}}
\newcommand{\RR}{\mathbb{R}}
\newcommandx*{\prob}[2][1,2={\omega}]{{{P^{#1}_{#2}}}}
\newcommandx*{\esper}[2][1,2={\omega}]{{{E^{#1}_{#2}}}}
\newcommand{\Prob}{\mathbb{P}}
\newcommand{\Esper}{\mathbb{E}}
\newcommand{\PROB}[1][]{{\mathbf{P}^{#1}}}
\newcommand{\ESPER}[1][]{{\mathbf{E}^{#1}}}

\newcommand{\elevation}{\mathfrak{E}}
\newcommand{\dirichlet}{\mathcal{E}}
\newcommand{\generator}{\mathcal{L}}
\newcommand{\depth}{\mathrm{depth}}
\newcommand{\domain}{\mathrm{Dom}}
\newcommand{\II}%
  {\mbox{$\mathrm{1\hspace{-3.1pt}l}$}} 
\newcommand{\argmax}{\mathop{\arg\max}}
\newcommand{\argmin}{\mathop{\arg\min}}
\newcommand{\locmin}{\mathop{\text{locmin}}}

\allowdisplaybreaks
\begin{document}


\title{Application of moderate deviation techniques to prove Sinai's Theorem on RWRE}
\author{Marcelo Ventura Freire\footnote{Corresponding address: mvf@usp.br;
the author was partially supported by FAPESP Grant 05/00248-6 and by CAPES}\\
Escola de Artes, Ci\^{e}ncias e Humanidades, \\Universidade de S\~{a}o Paulo (EACH/USP)
}

\maketitle

\begin{abstract}
We apply the techniques developed in \citet{Comets:Popov:2003} to present a new proof to Sinai's theorem \citep{Sinai:1982} on one-dimensional random walk in random environment (RWRE), working in a scale-free way to avoid rescaling arguments and splitting the proof in two independent parts: a quenched one, related to the measure $\prob$ conditioned on a fixed, typical realization $\omega$ of the environment, and an annealed one, related to the product measure $\Prob$ of the environment $\omega$.  
The quenched part still holds even if we use another measure (possibly dependent) for the environment.  

\noindent\textbf{Keywords:}
Random walk, random environment, Sinai's Walk, moderate deviations

\noindent\textbf{1991 MSC:} 60K37, 60G50
\end{abstract}


\section{Introduction}\label{sec:intro}
The Random Walk in Random Environment (RWRE) in $\ZZ$ is a jump process $\xi=\SET{\xi_t;t\in[0,\infty)}$ starting at $z\in\ZZ$ with law $\PROB[z]$ such that $\PROB[z](\cdot)=\int\prob[z](\cdot)\Prob(\dd\omega)$, where $\prob[z]$ is the law of a Markovian nearest-neighbor jump process starting at $z\in\ZZ$ with transition rates given by the fixed realization of the \emph{environment} $\omega=\SET{(\omega^-_x,\omega^+_x);x\in\ZZ}$, so that, for $h\searrow 0$, 
\begin{gather*}
  \prob[z](\xi_{t+h}=x\pm1|\xi_t=x)=\omega^\pm_x h+o(h),\\
  \prob[z](\xi_{t+h}=x|\xi_t=x)=1-(\omega^-_x+\omega^+_x)h+o(h),
\end{gather*}
and $\Prob$ is the law of the environment $\omega$, a product measure of the joint distribution of $\omega^-_0$ and $\omega^+_0$, so that the pairs $(\omega^-_x,\omega^-_x)$ are i.i.d.\ for $x\in\ZZ$.  
Expectations under $\PROB[x]$, $\Prob$, and $\prob[x]$ will be denoted as $\ESPER[x]$, $\Esper$, and $\esper[x]$ respectively and $\PROB[x]$ and $\prob[x]$ will be written $\PROB$ and $\prob$ when $x=0$.  

That model has been much studied in discrete time \citep[see][for an extensive review]{Zeitouni:RWRE} and recently in continuous time \citet{Comets:Popov:2003}, although  the discrete time model is embedded in the continuous time model, so there is no qualitative difference between them as long as the transition rates of the latter and the transition probabilities of the first are bounded away from $0$ and $\infty$ and from $0$ and $1$ respectively.  
A continuous state space version is introduced in \citet{Brox:1986} as the model of Brownian motion with random potential.  
Under $\PROB$, $\xi_t$ is not Markovian and the rates $\omega$ are homogeneous only at statistical level.  

\citet{Solomon:1975} established recurrence-transience criteria for the independent environment case, implying that $\xi_t$ is $\Prob$-a.s.\ recurrent if and only if 
$\Esper\log(\omega^+_0/\omega^-_0)=0$. 
Non-degenerate randomness of the environment is ensured if 
$0<\sigma^2:=\Esper\log^2(\omega^+_0/\omega^-_0)<\infty$,
so that RWRE is not a time-change of a simple random walk.  
Those conditions are called \emph{Sinai's regime}.  
The existence of a constant $\kappa>1$ such that
$\Prob(\kappa^{-1}\le\omega^\pm_0\le\kappa)=1$
is called \emph{ellipticity} and it is what ensures irreducibility of the RWRE and qualitative equivalence between the discrete and continuous time versions, besides finite variance for $\log\omega^+_0/\omega^-_0$ from Sinai's regime.  
Under these conditions, 
\begin{equation}\label{eq:cond}
  \Esper\log\frac{\omega^+_0}{\omega^-_0}=0, 
  \quad
  0<\sigma^2:=\Esper\log^2\frac{\omega^+_0}{\omega^-_0}<\infty,
  \quad
  \Prob(\kappa^{-1}\le\omega^\pm_0\le\kappa)=1,
\end{equation}
\citet{Sinai:1982} proved $\xi_t$ is of order $\log^2 t$, characterizing the strong sub-diffusive behavior of the RWRE in \eqref{eq:sinai}.  

\citet{Comets:Popov:2003} developed a new probabilistic approach which uses the KMT construction \citep{Komlos:Major:Tusnady:1975,Komlos:Major:Tusnady:1976} to study the moderate deviation of $\xi_t$ under $\PROB$ , but their techniques can be used to address questions such as extending Sinai's theorem to beyond environments with independent distribution.  

This new proof of Sinai's theorem separates in two parts what is due to the typical behavior of the random walk $\xi_t$ under $\prob$ for a fixed typical environment $\omega$ (the quenched part) and what is due to the typical behavior of the random environment $\omega$ under $\Prob$ (the annealed part).  
In the independent case, $\Prob$ is a product measure and the conjunction of Sinai's regime with ellipticity is sufficient condition to ensure that \eqref{eq:sinai} holds, but it is no longer sufficient in the dependent case.  
The quenched part of our proof is still valid in the dependent case, so that one needs to adapt only the annealed part for a dependent law for $\omega$ whose potential $V$ (defined ahead) still satisfy some suitable conditions.  
In this paper we present the proof for independent case and leave for a future paper the extension to dependent case.  

Another proof to Sinai's theorem has been given by \citet{Andreoletti:2005}, with a powerful approach, following the lines of \citet{Andreoletti:2006} and \citet{Andreoletti:2007}, 
where they strengthen the results of \citet{Sinai:1982} for the recurrent case still within Sinai's original conceptual framework, which included the creation of a hierarchy of refinements of valleys (or wheels) in the potential.  

Instead of investigating further the independent environment setup, our aim is to prove Sinai's theorem in a way we can extend the result to the case where $\Prob$ is no longer the product measure, like the recent extension of KMT construction to the dependent scenario in \citet{Berkes:Liu:Wu:2014} would allow, or when the potential converge to other stable L\'{e}vy processes than the Brownian motion.  

Our approach uses the fact that the potential converges weakly to a Brownian motion.  
We deal with the limiting Brownian motion coupled to the potential and then we are able to avoid rescaling arguments and work directly with the limit valleys in a scale-free fashion.  

In the next section we present the statement of Sinai's theorem; in Section~\ref{sec:notation} we define the concepts and notations we use; in Sections~\ref{sec:quenched:part} and \ref{sec:annealed:part} we give the proof, and in the appendix we present the proofs of the intermediate results needed in the Sections~\ref{sec:quenched:part} and \ref{sec:annealed:part}.  

\section{Main result}\label{sec:result}

Under Sinai's regime and ellipticity assumption we present an alternative proof of Sinai's theorem separated in two independent parts.  
In the quenched part we prove that a rescale of $\xi_t$ converges uniformly in $\prob$-probability as $t\to\infty$ to the same rescale of some process $m_t=m_t(\omega)$ function of the environment $\omega$ alone for any fixed typical environment $\omega$.  
In the annealed part we prove that the $\Prob$-measure of the set of typical environments $\omega$ converges to $1$.  

The Sinai's theorem can be rephrased as follows
\begin{theo}\label{theo:sinai}
If \eqref{eq:cond} holds, 
then there exists a jump process 
$
\{m_t; t\in[e,\infty)\}
$ 
such that, for any $\delta>\eps>0$, 
\begin{equation}\label{eq:quench:bound}
  \lim_{t\to\infty}\inf_{\omega\in\Typical}
    \prob\left(\left|\frac{\xi_t-m_t}{\log^2 t}\right|<\delta\right)
    =1,
\end{equation}
where $\Typical$ is 
such that 
\begin{equation}\label{eq:gamma:bound}
 \lim_{t\to\infty}\lim_{\eps\to 0}\Prob(\Typical)
   =\lim_{\eps\to 0}\lim_{t\to\infty}\Prob(\Typical)
   =1.  
\end{equation}
\end{theo}


The original formulation of the Sinai's theorem comes by as the following
\begin{coro}\label{coro:sinai}
  For any $\delta>0$,
  \begin{equation}\label{eq:sinai}
    \lim_{t\to\infty}
      \PROB\left(\left|\frac{\xi_t-m_t}{\log^2 t}\right|
                   >\delta\right)=0,
	\end{equation}
\end{coro}
immediately from 
\[
\PROB(|\xi_t-m_t|/\log^2 t>\delta) 
	\le\int_{\Typical} \prob(|\xi_t-m_t|/\log^2 t>\delta)\Prob(\dd\omega) +\Prob(\overline{\Typical}),
\] 
together with \eqref{eq:quench:bound} and \eqref{eq:gamma:bound}.  

\section{Notation and definitions}\label{sec:notation}

Transitions occur only between nearest neighbors, then the detailed balance equation $\theta_x\omega^+_x=\theta_{x+1}\omega^-_{x+1}$ can be solved, giving the \textbf{reversible measure} $\theta$
\[
  \theta_x=
  \begin{cases}
    \prod_{i=0}^{x-1}\frac{\omega^+_i}{\omega^-_{i+1}}, &x>0,\\
    1, &x=0\\
    \prod_{i=x}^{-1}\frac{\omega^-_{i+1}}{\omega^+_i}, &x<0
  \end{cases}
\]
that satisfies also $\theta_x\prob[x](\xi_t=y) =\theta_y\prob[y](\xi_t=x)$ for every $x,y\in\ZZ$ and $t>0$.  
Given a realization $\omega$, we define the \textbf{potential} $V=V[\omega]$ with domain $\ZZ$ as
\[
  V(x)=
  \begin{cases}
    \sum_{i=1}^{x}\log\frac{\omega^-_i}{\omega^+_i}, &x>0,\\
    0, &x=0\\
    \sum_{i=x+1}^{0}\log\frac{\omega^+_i}{\omega^-_i}, &x<0.  
  \end{cases}
\]
Ellipticity causes the rates to be bounded away from $0$ and $\infty$ and renders mutual domination between $\theta$ and $V$, for there exist positive constants $K_1,K_2$ such that $K_1 e^{-V(x)}\le \theta_x\le K_2 e^{-V(x)}$ for all $x$.  
Note that the function $w^{(n)}(t)$ of \citet{Sinai:1982} is our potential $V$ completed by linear interpolation and rescaled to converge weakly to a Brownian motion, so that $V(x)=w^{(n)}(x/\log^2 n)\log n$ for $x\in\ZZ$.

By hypothesis, the potential $V$ is a sum of i.i.d.r.v.'s with zero mean and finite second moment, therefore $V$ behaves like a random walk.  
By Donsker's Invariance Principle, $V(x\log^2 n)/\log n$ converges weakly as $n\to\infty$ to a two-sided Brownian motion $W(x)$ with diffusion coefficient $\sigma^2=\Esper(\log^2 \omega^-_0/\omega^+_0)$.  
We will use the strong approximation Theorem 1B of \citet[the KMT or hungarian construction]{Komlos:Major:Tusnady:1976} to work directly with the limiting Brownian motion $W$ (which possesses the self-scaling property) in substitution of the potential $V$.  

Accordingly, in a possibly enlarged probability space there exist a version of our environment process $\omega$ and a two-sided Brownian motion $W$ with diffusion constant $\sigma$ such that
for some $\hat\kappa>0$
\begin{equation}\label{eq:KMT}
  \Prob\left(\limsup_{x\to\pm\infty}
             \frac{|V(x)-W(x)|}{\log|x|}\le\hat\kappa\right)=1.  
\end{equation}

\citet{Sinai:1982} worked the idea of \emph{refinement} of the function $w^{(n)}(t)$ while we will work the idea introduced by \citet{Comets:Popov:2003} of \emph{$t$-stable wells} and \emph{$t$-stable points} on the potential $V$ and on its scaling limit $W$.  

We can define the concept of \emph{$t$-stability} for any real function $f$ with domain $\domain(f)$ (which may be either $V$ or $W$ with domains $\ZZ$ or $\RR$ resp.), but we need first some previous definitions.  
In the following definitions and whenever necessary, we consider all maxima, minima, suprema and infima of $f$ over a set $I$ as over $I\cap\domain(f)$.

\begin{figure}
  \centering
  \resizebox{\textwidth}{!}{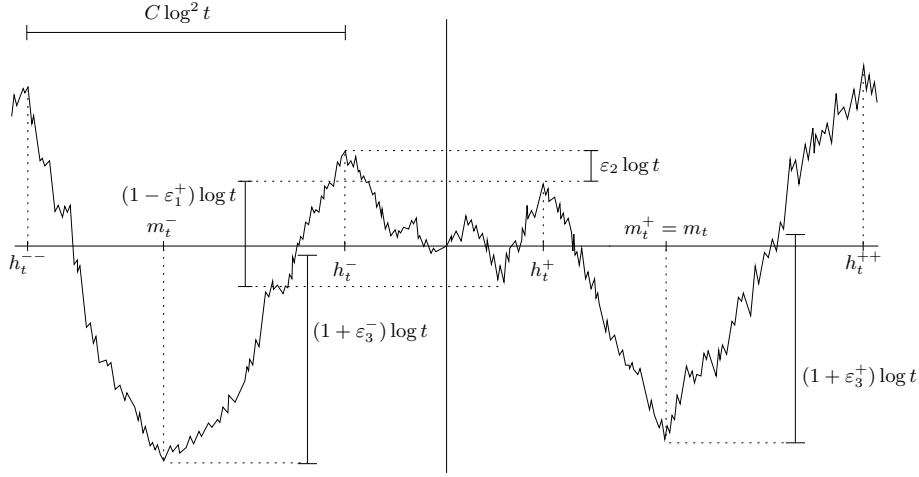}
  \caption{A function $f$ with two $t$-stable wells}
  \label{fig:t:stable:wells}
\end{figure}

We say that a finite interval $I=[a,b]$ is a \textbf{well} on a function $f$ if $a=\argmax_{x\in[a,c]}f(x)$ and $b=\argmax_{x\in[c,b]}f(x)$, where $c=\argmin_{x\in[a,b]}f(x)$.  
We define the \textbf{depth} of a well $I=[a,b]$ on $f$ as $\depth(I):=\min\{f(a),f(b)\} -\min_{x\in[a,b]}f(x)$.  
For $t>1$, we say that a point $m\in\domain(f)$ is a \textbf{$t$-stable point} of $f$ if $m=\argmin_{x\in[l,r]}f(x)$, where $l=l(t,m):=\sup\{x\in(-\infty,m];f(x)\ge f(m)+\log t\}$ and $r=r(t,m):=\inf\{x\in[m,\infty);f(x)\ge f(m)+\log t\}$.  
In plain words, a $t$-stable point is the bottom of a well at least as deep as $\log t$, as the points $m^-_t$ and $m^+_t$ in Figure~\ref{fig:t:stable:wells}.  

We define also the set $\Stable$ of all $t$-stable points of $f$ and let $\Stable^+:=\Stable\cap[0,\infty)$ and $\Stable^-:=\Stable\cap(-\infty,0]$.  
For any $t>e$, $\Stable{[V]}$, $\Stable{[W]}$, and their traces $\Stable\cap(-\infty,x)$ and $\Stable\cap(x,\infty)$ 
are infinite.  
Besides, all their elements are isolated points both for $V$ and for $W$, because in one case $\domain(V)$ is an isolated point set and in the other case, between its local minima, $W$ need to raise and fall both at least $\log t$ before another local minimum can belong to $\Stable{[W]}$, so an accumulation point in $\Stable{[W]}$ $\Prob$-a.s.\ cannot occur.  

Between two successive $t$-stable points $m$ and $m'$, there exists a peak $h=\argmax_{x\in[m,m']}f(x)$ separating two adjacent well of depth of at least $\log t$, so we define the set $\Peaks$ of peaks of $f$ which separate $t$-stable points as $\Peaks:=\{h\in\domain(f); \exists\, m,m'\in\Stable: h=\argmax_{x\in[m,m']}f(x)\}$.  

We define the \textbf{$t$-stable well} $\Well$ of the $t$-stable point $m\in\Stable$ as $\Well(m):=[\max\Peaks\cap(-\infty,m), \min\Peaks\cap(m,\infty)]$, so that any $t$-stable well is formed by two successive $h,h'\in\Peaks$ with only one $m\in\Stable$ in between.  


For the proofs, we define the $t$-stable points which are closest to the origin as well as the peaks surrounding them as
\begin{align*}
 \begin{gathered}
  m^-_t    := \max\Stable^- \\
  h^-_t    := \argmax_{x\in[m^-_t,0]}f(x) \\
  m^{--}_t := \max\Stable\cap(-\infty,h^-_t) \\
  h^{--}_t := \argmax_{x\in[m^{--}_t,m^-_t]}f(x)
 \end{gathered}
 \quad\quad\quad\quad
 \begin{gathered}
  m^+_t    := \min\Stable^+\\
  h^+_t    := \argmax_{x\in[0,m^+_t]}f(x)\\
  m^{++}_t := \min\Stable\cap(h^+_t,\infty)\\
  h^{++}_t := \argmax_{x\in[m^+_t,m^{++}_t]}f(x).  
 \end{gathered}
\end{align*}
These definitions are illustrated in Figure~\ref{fig:t:stable:wells} for $f=W$, but notice that, while $h^-_t\in\Peaks$ in that example, we also have that $h^+_t\not\in\Peaks$, for $f(h^+_t)$ is not the maximum between $m^-_t$ and $m^+_t$.  Still in Figure~\ref{fig:t:stable:wells} we can see the $t$-stable wells $\Well(m^-_t)=[h^{--}_t,h^-_t]$ and $\Well(m^+_t)=[h^-_t,h^{++}_t]$.

Now, we can define the jump process $\{m_t; t>e\}$ from \eqref{eq:quench:bound} that will attract the random walk $\xi_t$ at each moment $t>e$ as 
\[
  m_t:=
  \begin{cases}
    m^-_t, &\text{if $f(h^+_t)>f(h^-_t)$}\\
    m^+_t, &\text{if $f(h^+_t)<f(h^-_t)$},
  \end{cases}
\]
i.e., $m_t$ will be the closest-to-the-origin $t$-stable point for each instant $t$.  

\begin{figure}
  \centering
  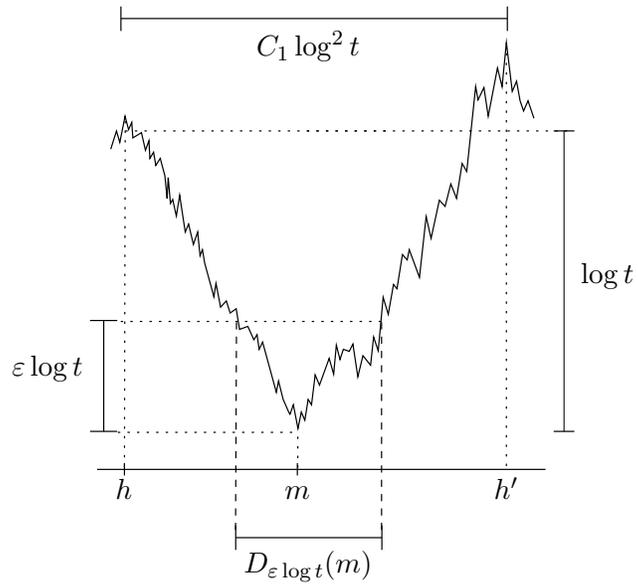
  \caption{At the bottom of a $t$-stable well}
  \label{fig:bottom:well}
\end{figure}

For $m\in\Stable$ with $\Well(m)=[h,h']$, $h,h'\in\Peaks$, and $0<a\le\depth(\Well(m))$, we define the \textbf{$a$-neighborhood} $\Neigh[a]$ of $m$ as
$
  \Neigh[a](m):=[\mathfrak{l}(m,a),\mathfrak{r}(m,a)],
$
where we have $\mathfrak{l}(m,a):=\inf\SET{x\in[h,m]: f(x)-f(m)<a}$ and $\mathfrak{r}(m,a):=\sup\SET{x\in[m,h']: f(x)-f(m)<a}$.  
Notice that $\Neigh[a](m)\subset\Well(m)$, since $h\le\mathfrak{l}<\mathfrak{r}\le h'$ by definition,
and also that $f(x)-f(m)>\eps\log t$ for $x\in\Well(m)\smallsetminus\Neigh(m)$.  
An instance of an $(\eps\log t)$-neighborhood $\Neigh(m)$ of a $t$-stable point $m$ is shown in Figure~\ref{fig:bottom:well}.  

\begin{figure}
 \centering
 \resizebox{.8\textwidth}{!}{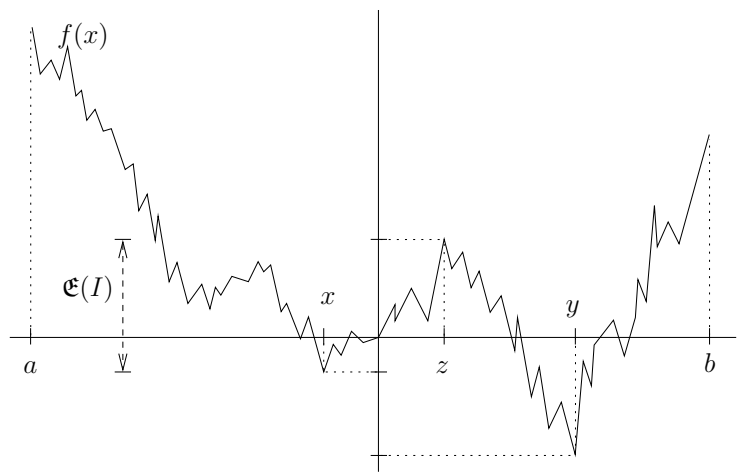}
 \caption{Elevation $\elevation(I)$ of a function $f$ 
          over the interval $I=[a,b]$}
 \label{fig:elevation}
\end{figure}

\citet{Mathieu:1994} 
defined the \textbf{elevation} 
$\elevation$ 
of $f$ in the interval $I=[a,b]$ as 
\[
\elevation(I):=\max_{x,y\in I} \max_{z\in\mathcal{I}(x,y)}(f(z)-f(x)-f(y))+\min_{v\in I}f(v)
\]
or, equivalently in our case, 
\[
\elevation(I)=\max_{x\in \locmin(f,I)}\max_{z\in\mathcal{I}(x,x_{\min})}(f(z)-f(x)),
\]
where $\mathcal{I}(a,b)=[a,b]\cup[b,a]$, $x_{\min}=\argmin_{v\in I} f(v)$ is the global minimum of $f$ over $I$ and $\locmin(f,I)$ is the set of local minima of $f$ over $I$ except the global minimum $x_{\min}$.  
For $I\subset J$, we have $\elevation(I) \le\elevation(J)$.  
The definition is illustrated in Figure~\ref{fig:elevation}.  

The above definitions can be used with both the actual potential $V$ of the environment and its scaling limit $W$, the Brownian motion coupled to $V$ accordingly to \eqref{eq:KMT}.  
Since we intend to use $W$ in the place of $V$, all $t$-stable points, $t$-stable wells and $a$-neighborhoods from this point on will be relative to the Brownian motion $W$ unless explicitly stated in notation.  

We must draw attention to the fact that the points $m^\pm_t$, $h^\pm_t$, $m^{\pm\pm}_t$, and $h^{\pm\pm}_t$ relative to $W$ defined above are $\Prob$-a.s.\ non integers. 
So, throughout this paper, statements like ``the random walk $\xi$ hits a $t$-stable point $m$'' 
means that it hits the site $x\in\ZZ$ which is closest to $m$.  
Throughout this paper, real points $x$ will be replaced, if the context requires, with the closest integer, so that we may still denote by the same symbol $x$, if no confusion can occur.

At last, $K_1, K_2,...$ denote positive constants that may change from line to line.  

\section{Quenched part of the proof}\label{sec:quenched:part}

Putting aside technicalities, the idea of this part of proof is that, for any typical environment $\omega$, this is what  happens with large probability:
$(i)$ the particle will leave the interval $[h^-_t,h^+_t]$ before the instant $t$; 
$(ii)$ the particle will choose to leave $[h^-_t,h^+_t]$ through the lowest of the peaks in direction of its respective $t$-stable point; 
$(iii)$ prior to instant $t$, the particle will reach the $t$-stable point, that will be either $m^-_t$ or $m^+_t$ depending on the lowest of $W(h^-_t)$ and $W(h^+_t)$; 
$(iv)$ once the $t$-stable point is reached before $t$, the particle will not leave the $t$-stable well until the instant $t$; 
$(v)$ still within the $t$-stable well until the instant $t$, the particle will oscillate inside a rather narrow neighborhood of the $t$-stable point; 
$(vi)$ the breadth of that neighborhood scaled by $\log^2 t$ will be arbitrarily small for $t$ large enough.  


Fix $M>2$ arbitrarily and consider $t>e$. 
Using the Brownian motion $W$ coupled to the potential $V$ of \eqref{eq:KMT}, 
let $\Typical[1][t]$ be the set of environments $\omega$ whose potential $V$ is close enough to $W$ within the radius $|x|\le\log^M t$ 
and let $\Typical[2][t]$ be the set of environments $\omega$ whose two $t$-stable wells surrounding the origin are within the radius $|x|\le\log^M t$
\begin{gather} 
  \Typical[1][t]:=
    \SET{|V(x)-W(x)|<\hat\kappa M\log\log t,|x|<\log^M t}, \label{eq:gamma:KMT} \\
  \Typical[2][t]:=
    \SET{|h_t^{--}|<\log^M t, |h_t^{++}|<\log^M t}, \label{eq:gamma:contained}
\end{gather}
where $\hat\kappa$ in \eqref{eq:gamma:KMT} comes from \eqref{eq:KMT}.  
Here, $W$ is the Brownian motion coupled with the potential $V$ through KMT construction \eqref{eq:KMT}, so that we are able to use either $V$ or $W$, whichever is easier to deal with in context.  


Let $\tau_A:=\inf\{t>0:\xi_t\in A\}$ be the \emph{hitting time} of $\xi$ in $A\subset\ZZ$ (with $\tau_x=\tau_{\{x\}}$) and consider the events
\begin{align*}
  A_1&:=\SET{\tau_{\{m^-_t,m^+_t\}}<t}
  &   A^\pm_3&:=\SET{\tau_{H^\pm_t}>t}\\
  A^\pm_2&:=\SET{\tau_{\{m^-_t,m^+_t\}}=\tau_{m^\pm_t}}
  &   A^\pm_4&:=\SET{\xi_t\in\Neigh(m^\pm_t)}.  
\end{align*}
where $m^\pm_t$, $h^\pm_t$, $h^{\pm\pm}_t$, and $\Neigh[a](m^\pm_t)$ are relative to the scaling limit $W$ of the potential $V$ and $H^+_t:=\{h^-_t,h^{++}_t\}$ and $H^-_t:=\{h^{--}_t,h^+_t\}$ are the peaks around $m^\pm_t$ respectively.

Then we have
\begin{multline}\label{eq:quenched:bound:1}
    \prob\left(|\xi_t-m^\pm_t|<|\Neigh(m^\pm)|\right)
      \ge\prob(A^\pm_4)
      \ge\prob(A_1,A^\pm_2,A^\pm_3,A^\pm_4)\\
      \ge 1 - \prob\left(\overline{A_1}\right) 
            - \prob\left(\overline{A^\pm_2}\right)
            - \prob\left(\overline{A^\pm_3}\middle|A_1,A^\pm_2\right)
            - \prob\left(\overline{A^\pm_4}\middle|A_1,A^\pm_2,A^\pm_3\right).  
\end{multline}
Such probabilities can be bounded through the next four lemmas, whose proofs are left for the appendices of this paper.  

\begin{lemma}\label{lemma:leave:subwell}
For $\omega\in\Typical[1][t]\cap\Typical[2][t]$ and $t$ large enough,
\begin{equation}\label{eq:bound:term1}
  \prob\left(\tau_{\{m^-_t,m^+_t\}}>t\right) 
    \le t^{-\eps^+_1}+t^{-\eps^-_1}
    \le K_1 t^{-\eps_1},
\end{equation}
where $\eps^\pm_1=1-\elevation(\Well(m^\pm_t))/\log t$ and $\eps_1=\min\{\eps^-_1,\eps^+_1\}$.  
\end{lemma}

\begin{lemma}\label{lemma:choose:leave}
For $\omega\in\Typical[1][t]\cap\Typical[2][t]$,
if $W(h^-_t)\lessgtr W(h^+_t)$, 
then 
  \begin{equation}\label{eq:bound:term2}
    \prob\left(\tau_{\{m^-_t,m^+_t\}}=\tau_{m^\pm_t}\right)
      \le t^{-\eps_2} \log^{(2\hat\kappa+1)M} t,
  \end{equation}
where $\eps_2=|W(h^-_t)-W(h^+_t)|/\log t$.  
\end{lemma}

We state that 
\begin{equation}\label{eq:bound:term3}
  \prob\left(\overline{A^\pm_3}\middle|A_1,A^\pm_2\right) 
    \le\prob[m^\pm_t]\left(\overline{A^\pm_3}\right)
    \le K_3 t^{-\eps^\pm_3}\log^{2\hat\kappa M} t,
\end{equation}
where 
$\eps^\pm_3={\depth(\Well(m^\pm_t))}/{\log t}-1$, 
because
\begin{align*}
  \prob\left(\overline{A^\pm_3},A_1,A^\pm_2\right)
    &=\int_{[0,t)}
       \prob\left(\tau_{H^\pm_t}<t \middle|
             \tau_{m^\pm_t}=s,A^\pm_2\right)
       \dd\prob\left(\tau_{m^\pm_t}<s,A^\pm_2\right)\\
    &=\int_{[0,t)}
       \prob[m^\pm_t]\left(\tau_{H^\pm_t}<t-s\right)
       \dd\prob\left(\tau_{m^\pm_t}<s,A^\pm_2\right)\\
    &\le\int_{[0,t)}
       \prob[m^\pm_t]\left(\tau_{H^\pm_t}<t\right)
       \dd\prob\left(\tau_{m^\pm_t}<s,A^\pm_2\right)\\
    &\le\prob[m^\pm_t]\left(\overline{A^\pm_3}\right)
       \cdot\prob\left(\tau_{m^\pm_t}<t,A^\pm_2\right)
     =\prob[m^\pm_t]\left(\overline{A^\pm_3}\right)\cdot\prob\left(A_1,A^\pm_2\right)
\end{align*}
and 
\begin{lemma}\label{lemma:stay:well}
For $\omega\in\Typical[1][t]\cap\Typical[2][t]$, $m\in\Stable$, and $\Well(m)=[h,h']$ with $h,h'\in\Peaks$, we have that
  \begin{equation*}
    \prob[m]\left(\tau_{\{h,h'\}}<t\right)
      \le K_3\cdot t^{-\eps_3}\cdot\log^{2\hat\kappa M} t
  \end{equation*}
where $\eps_3={\depth(\Well(m))}/{\log t}-1$.  
\end{lemma}

Finally, we also state that
\begin{equation}\label{eq:bound:term4}
  \prob\left(\overline{A^\pm_4}\middle|A_1,A^\pm_2,A^\pm_3\right)
    \le K_4 \; t^{-\eps} \log^{(2\hat\kappa+1)M} t.  
\end{equation}
because
\begin{align*}
  \prob&\left(\overline{A^\pm_4},A_1,A^\pm_2,A^\pm_3\right)
     =\prob\left(\xi_t\not\in\Neigh(m^\pm_t),
         \tau_{m^\pm_t}< t,A^\pm_2,A^\pm_3\right)\\
    &=\int_{[0,t)}
       \prob\left(\xi_t\not\in\Neigh(m^\pm_t)\middle|
             \tau_{m^\pm_t}=s,A^\pm_2,A^\pm_3\right)
       \dd\prob\left(\tau_{m^\pm_t}<s,A^\pm_2,A^\pm_3\right)\\
    &=\int_{[0,t)}
        \prob[m^\pm_t]\left(\xi_{t-s}\not\in\Neigh(m^\pm_t)\middle|
                     \tau_{H^\pm_t}>t-s\right)
        \dd\prob\left(\tau_{m^\pm_t}<s,A^\pm_2,A^\pm_3\right)\\
    &\le \int_{[0,t)}
        K_1 \; t^{-\eps} \log^{(2\hat\kappa+1)M} t \;
        \dd\prob\left(\tau_{m^\pm_t}<s,A^\pm_2,A^\pm_3\right)\\
    &=K_1 \; t^{-\eps}\log^M t \cdot 
      \prob(A_1,A^\pm_2,A^\pm_3),
\end{align*}
with the inequality due to
\begin{lemma}\label{lemma:oscillating:bottom}
For $\omega\in\Typical[1][t]\cap\Typical[2][t]$, 
if $m\in\Stable$ and $\Well(m)=[h,h']$, $h,h'\in\Peaks$, 
then for $s<t$
\[
 \prob[m]\left(\xi_s\not\in\Neigh(m)\middle|\tau_{\{h,h'\}}>s\right)
   \le K_1 \; t^{-\eps} \log^{(2\hat\kappa+1)M} t.  
\]
\end{lemma}

Gathering \eqref{eq:bound:term1}--\eqref{eq:bound:term4} and applying them into \eqref{eq:quenched:bound:1} gives
\begin{multline}\label{eq:quenched:bound:2}
  \prob\left(\left|\xi_t-m^\pm_t\right|<\left|\Neigh(m^\pm_t)\right|\right)\\
    \ge 1-K_1 \; t^{-\eps_1}
         - t^{-\eps_2} \log^{(2\hat\kappa+1)M} t,
         -K_2 \; t^{-\eps^\pm_3}\log^{2\hat\kappa M} t
         -K_3 \; t^{-\eps} \log^{(2\hat\kappa+1)M} t.
\end{multline}

In order to control the loose terms $\eps_1$, $\eps_2$, and $\eps_3$ above, 
we define the set 
$\Typical[3]$ of all $\omega$ whose difference between the height of first peaks around the origin is large enough,
$\Typical[4,\pm]$ of all $\omega$ whose $t$-stable well's elevation is smaller enough than $\log t$, and 
$\Typical[5,\pm]$ of all $\omega$ whose $t$-stable well's depth is larger enough than $\log t$:
for $t>e$ and $\eps\in(0,\delta)$
\begin{gather}
  \Typical[3]:=
    \SET{\frac{|W(h_t^-)-W(h_t^+)|}{\log t}>\eps}
         \label{eq:gamma:peak:peak}\\
  \Typical[4,\pm]:=
    \SET{\frac{\log t - \elevation(\Well(m^\pm_t))}{\log t}>\eps}
         \label{eq:gamma:elevation}\\
  \Typical[5,\pm]:=
    \SET{\frac{\depth(\Well(m^\pm_t))-\log t}{\log t}>\eps}.
         \label{eq:gamma:depth}
\end{gather}

Now we use the fact that, for $\omega\in\Typical[1][t]\cap\Typical[2][t]\cap\Typical[3]\cap\Typical[4,\pm]\cap\Typical[5,\pm]$, \eqref{eq:quenched:bound:2} reduces to
\begin{equation}\label{eq:quenched:bound:3}
  \prob\left(\left|\xi_t-m^\pm_t\right|<\left|\Neigh(m^\pm_t)\right|\right)
    \ge 1 - K_4 \; t^{-\eps} (1+\log^{2\hat\kappa M} t+2\log^{(2\hat\kappa+1)M} t).
\end{equation}

To control the breadth of $\Neigh(m^\pm_t)$, we consider also, for $t>e$ and $\eps\in(0,\delta)$
\begin{gather}
  \Typical[6,\pm]:=
    \SET{|\Neigh(m^\pm_t)|<\eps\log^2 t}.
         \label{eq:gamma:neighborhood}
\end{gather}

Once we have $|\xi_t-m^\pm_t|<|\Neigh(m^\pm_t)|<\eps\log^2 t$
for $\omega\in\Typical[6,\pm]$, we also have
\begin{multline*}
\prob\left(|\xi_t-m^\pm_t|<|\Neigh(m^\pm_t)|\right)
  \le \prob\left(|\xi_t-m^\pm_t|<\eps\log^2 t\right)\\
  \le \prob\left(\frac{|\xi_t-m^\pm_t|}{\log^2 t}<\eps\right)
  \le \prob\left(\frac{|\xi_t-m^\pm_t|}{\log^2 t}<\delta\right).
\end{multline*}

Finally, applying \eqref{eq:quenched:bound:3} in the inequality above, we have~\eqref{eq:quench:bound}, since
\begin{equation*}
\prob\left(\frac{|\xi_t-m^\pm_t|}{\log^2 t}<\delta\right)
  \ge 1 - K_4 \; t^{-\eps} (1+\log^{2\hat\kappa M} t+2\log^{(2\hat\kappa+1)M} t),
\end{equation*}
which converges to 1 as $t\to\infty$ for any $\eps\in(0,\delta)$ 
for any $\omega\in\Typical$ given by
\begin{equation}\label{def:gamma}
  \Typical:=
    \Typical[1][t]
    \cap\Typical[2][t]
    \cap\Typical[3]
    \cap\Typical[4,-]\cap\Typical[5,-]\cap\Typical[6,-]
    \cap\Typical[4,+]\cap\Typical[5,+]\cap\Typical[6,+].
\end{equation}

\section{Annealed part of the proof}\label{sec:annealed:part}

Now we prove that the $\Prob$-measure of every set in \eqref{def:gamma} above converges to $1$, so that $\Prob(\Typical)\to1$ as $t\to\infty$ and $\eps\to0$.  

To prove the convergence of $\Prob(\Typical[1][t])$, we notice that \eqref{eq:KMT} assures that
\[
\Prob(\Typical[1][t])
=
\Prob\left(\max_{x\in[-\log^M t,\log^M t]} \frac{|V(x)-W(x)|}{M\log\log t}\le\hat\kappa\right)
\cvg{t\to\infty}1.
\]

%
%
%
%

To prove the convergence for $\Typical[2][t]$ to $\Typical[5,\pm]$, we use this 
\begin{propo}\label{propo:Brownian:scaling}
Let $W$ be a Brownian motion and $\hat W(\cdot)=a\,W(\cdot/a^2)$ be a rescaling of $W$.  
Then, for $a,b>0$, $t>e$ and $m\in\Stable(W)$
\[
  \Stable[t^a](\hat W)=a^2\Stable[t](W);
  \quad
  h^\pm_{t^a}(\hat W)=a^2 h^\pm_t(W);
  \quad
  \Neigh[ab](a^2m)(\hat W)=a^2\Neigh[b](m)(W).  
\]
\end{propo}
The proof is immediate from definitions and standard scaling arguments, so it is omitted.  
As an immediate consequence,  $a^2\Peaks[t](W)=\Peaks[t^a](\hat W)$, since $W\eqdistr-W$ renders $\Peaks(W)\eqdistr\Stable(-W)$.  

Applying Proposition~\ref{propo:Brownian:scaling} above with $a=1/\log t$ gives 
\begin{align}
h^{\pm\pm}_t/\log^2 t &\eqdistr h^{\pm\pm}_e,
  \label{eq:distrib:well:breadth}\\
W(h^\pm_t)/\log t &\eqdistr W(h^\pm_e),  
  \label{eq:distrib:peak:height}\\
\elevation(\Well(m^\pm_t))/\log t &\eqdistr\elevation(\Well[e](m^\pm_e)), 
  \label{eq:distrib:well:elevation}\\
\depth(\Well(m^\pm_t))/\log t &\eqdistr\depth(\Well[e](m^\pm_e)),
  \label{eq:distrib:well:depth}\\
\Neigh(m^\pm_t)/\log^2 t &\eqdistr\Neigh[\eps](m^\pm_e),
  \label{eq:distrib:neighborhood:breadth}
\end{align}
whose right-hand-side's distributions do not depend on $t$ and, except in~\eqref{eq:distrib:neighborhood:breadth}, do not depend on $\eps$.

According to~\eqref{eq:distrib:well:breadth},
\begin{multline*}
\Prob(\Typical[2][t])
  \ge 1 - \Prob\left(|h^{--}_t|\ge\log^M t\right) - \Prob\left(|h^{++}_t|\ge\log^M t\right)\\
  \ge 1 - 2\Prob\left(\frac{|h^{++}_t|}{\log^2 t}\ge\log^{M-2} t\right)
  \ge 1 - 2\Prob\left(|h^{++}_e|\ge\log^{M-2} t\right),
\end{multline*}
which does not depend on $\eps$ and converges to 1 as $t\to\infty$.

According to~\eqref{eq:distrib:peak:height}--\eqref{eq:distrib:well:depth}, the distribution of the fractions in \eqref{eq:gamma:peak:peak}--\eqref{eq:gamma:depth} depend only on $\eps$, so that the probabilities $\Prob(\Typical[3])$, $\Prob(\Typical[4,\pm])$, and $\Prob(\Typical[5,\pm])$ also depend only on $\eps$.  Since the fractions inside \eqref{eq:gamma:peak:peak}, \eqref{eq:gamma:elevation} and \eqref{eq:gamma:depth} are strictly positive r.v.'s with absolute continuous distributions, then $\Prob(\Typical[3])$, $\Prob(\Typical[4,\pm])$ and $\Prob(\Typical[5,\pm])$ converge to $1$ as $\eps\to0$.  

Since $|\Neigh[\eps](m^\pm_e)|=O(\eps^2)=o(\eps)$ by scaling properties of Brownian motion,  \eqref{eq:distrib:neighborhood:breadth} gives
$
  \Prob(\Typical[6,\pm])
    =\Prob(\Neigh(m^\pm_t)/\log^2 t<\eps) 
    =\Prob(|\Neigh[\eps](m^\pm_e)|<\eps), 
$
which does not depend on $t$ and converges to 1 as $\eps\to0$.

Finally we get \eqref{eq:gamma:bound} for $\Prob(\Typical[1][t])$ and $\Prob(\Typical[2][t])$ 
are constant in relation to $\eps$ and 
converge to $1$ as $t\to\infty$ and 
$\Prob(\Typical[3])$, $\Prob(\Typical[4,\pm])$, $\Prob(\Typical[5,\pm])$, and $\Prob(\Typical[6,\pm])$ 
are constant in relation to $t$ and 
converge to 1 as $\eps\to0$.  


\bibliographystyle{plainnat}
\bibliography{biblio}

\appendix
\section{Auxiliary results}\label{sec:aux:results}

\subsection{Reflected RWRE in an interval}\label{sec:reflected:rw}

In order to use the reversible measure $\theta$ of the RWRE $\xi$ in the proofs of following sections, we construct a version $\hat\xi_t$ of $\xi_t$ reflected in an interval $[a,b]$ and started at the same origin $y\in(a,b)$, through this following coupling.  

Let $\{U_n;n\in\NN^*\}$ and $\{V_n;n\in\NN^*\}$ be two independent sequences of i.i.d.r.v.'s with $\mathrm{Unif}(0,1)$ and $\mathrm{Expon}(1)$ distributions respectively.  
We define the process $\xi=\{\xi_t;t\in\RR^+\}$ and its sequence $\{T_n;n\in\NN\}$ of transition times by
\begin{gather*}
  \begin{alignedat}{2}
    \xi_0&:=y,& \qquad\qquad T_0&:=0, \\
    \xi_s&:=\xi_{T_{n-1}},\;\forall s<T_n,
             & T_n&:=T_{n-1}+V_n/(\omega^-_{\xi_{T_{n-1}}}+\omega^+_{\xi_{T_{n-1}}}),
  \end{alignedat}\\
  \xi_{T_n}:=\xi_{T_{n-1}}
      -\II\left(
        U_n<\frac{\omega^-_{\xi_{T_{n-1}}}}
                 {\omega^-_{\xi_{T_{n-1}}}
                  +\omega^+_{\xi_{T_{n-1}}}}
      \right)
     +\II\left(
        U_n>\frac{\omega^-_{\xi_{T_{n-1}}}}
                 {\omega^-_{\xi_{T_{n-1}}}
                  +\omega^+_{\xi_{T_{n-1}}}}
      \right)
\end{gather*}
and we define $\hat \xi_t$ and $\hat T_n$ analogously with the same $U_n$'s and $V_n$'s but with $\hat\omega$ instead of $\omega$, where $\hat\omega$ is such that $\hat\omega^\pm_x=\omega^\pm_x$ for $x\in(a,b)$ and reflected at the extremes $a,b$ with $\hat\omega^-_a=0$, $\hat\omega^+_a=\omega^+_a$, $\hat\omega^-_b=\omega^-_b$ and $\hat\omega^+_b=0$ and with $\hat\omega^\pm_x$ arbitrary for $x$ outside $[a,b]$.  
Let $\hat\tau_A:=\inf\{t>0:\hat\xi_t\in A\}$ the hitting time of $\hat\xi$, just as $\tau_A$ is the hitting time of $\xi$.  
In that construction, we can easily see that $\tau_{\{a,b\}}=\hat\tau_{\{a,b\}}$ and $\hat\xi_t=\xi_t$ for $t\le\tau_{\{a,b\}}$.  

The solution to the detailed balance equation for $\hat\xi$ is $\Prob$-a.s.\ summable, 
so $\hat\xi$ is $\Prob$-a.s.\ $\prob[][\hat\omega]$-positive-recurrent 
and we can find that the $\prob[][\hat\omega]$-stationary distribution $\mu=\mu_{[a,b]}$ of $\hat\xi$ is $\Prob$-a.s.\ $\mu(A)=\sum_{x\in A\cap[a,b]}\theta_x/\sum_{z\in[a,b]}\theta_z$.  
The potential $\hat V$ for $\hat\xi$ is $\hat V(x)=V(x)-V(y)$ for $x\in[a,b]$ and arbitrary outside $[a,b]$.  
As $\hat\xi$ is $\Prob$-a.s.\ $\prob[][\hat\omega]$-reversible, we have the symmetry of the \emph{infinitesimal generator} $\generator=\generator([a,b])$ of $\hat\xi$ given by
\[
  \generator f(x)
    :=\lim_{t\to 0}\frac{\esper[x][\hat\omega] f(\hat\xi_t)-f(x)}{t}
    =(f(x+1)-f(x))\cdot\hat\omega^+_x+(f(x-1)-f(x))\cdot\hat\omega^-_x
\]
and then we can define the \emph{Dirichlet form} $\dirichlet=\dirichlet([a,b])$ of $\hat\xi$ as $\dirichlet(f,f):=-\langle\generator f,f\rangle_{L^2(\mu)} =\sum_{x\in[a,b)}(f(x+1)-f(x))^2\hat\omega^+_x\mu(x)$ for any $f\in L^2(\mu)$ and the \emph{spectral gap} $\lambda=\lambda([a,b])$ of $\hat\xi$ as 
\begin{equation}\label{def:spectral:gap}
  \lambda 
    :=\inf\{\dirichlet(f,f): f\in L_2(\mu), 
                             \esper[][\mu] f(\hat\xi_0)=0, 
                             \esper[][\mu] f(\hat\xi_0)^2=1\}.  
\end{equation}
We can approximate the spectral gap $\lambda([a,b])$ with the elevation $\elevation[\hat V]([a,b])= \elevation[V]([a,b])$ of $\hat V$ over $[a,b]$ through Proposition 3.1 of \citet{Comets:Popov:2003} or II.0 of \citet{Mathieu:1994}: for $M>0$,
\begin{equation}\label{eq:spectral:gap:elevation}
  \lim_{t\to\infty}\sup_{I\subset[-\log^M t,\log^M t]}
    \frac{|\log\lambda(I)+\elevation(I)|}{\log t}
    = 0.  
\end{equation}

\subsection{Proof of Lemma~\ref{lemma:leave:subwell}}\label{sec:proof:lemma:subwell}

This an application of Lemma 3.1 in \citet{Comets:Popov:2003}, whose proof deals with the reflected version $\hat\xi$ of the RWRE $\xi$ introduced above.  
In adapted notation, it states that, for $\omega\in\Typical[1][t]$ and for every $x$ such that $m<x<m'$, for any two consecutive $t$-stable points $m,m'\in\Stable$ with the peak $h\in\Peaks$ in between, we have that
\begin{align*}
  &\prob[x](\tau_{\{m,m'\}}>t/y)\\
    &\le\exp\left\{
         -t^{\frac12\left(1-\frac{\elevation(I^+)}{\log t}\right)}
           \left(
             \frac{K_1}{\Delta\log^{2\hat\kappa}t}
             -K_2 e^{\gamma/2}
              \exp\left\{
                -\lambda(I^+)e^{\elevation(I^+)}
                 t^{\frac12\left(1-\frac{\elevation(I^+)}{\log t}\right)}/2y
              \right\}
           \right)
       \right\}\\
    &+\exp\left\{
         -t^{\frac12\left(1-\frac{\elevation(I^-)}{\log t}\right)}
           \left(
             \frac{K_1}{\Delta\log^{2\hat\kappa}t}
             -K_2 e^{\gamma/2}
              \exp\left\{
                -\lambda(I^-)e^{\elevation(I^-)}
                 t^{\frac12\left(1-\frac{\elevation(I^-)}{\log t}\right)}/2y
              \right\}
           \right)
       \right\}
\end{align*}
where $I^+=[h,m']$ and $I^-=[m,h]$, $\Delta=m'-m$, $\gamma=\max_{x\in[m,m']}V(x)-\min_{x\in[m,m']}V(x)$, $\lambda$ is the spectral gap introduced in \eqref{def:spectral:gap} and the constants $K_1$ and $K_2$ depend only on $\omega$.  

Here we will take $y=1$, $m=m^-_t$, $m'=m^+_t$, and $x=0$, take also as $h$ the only element of $\{h^-_t,h^+_t\}\cap\Peaks$ and consider that here $\omega\in\Typical[1][t]\cap\Typical[2][t]$,
which makes $\Delta\le 2\log^M t$ and makes $\gamma=V(h)-\min_{x=m^\pm_t}V(x)>\log t$, which makes $e^{\gamma/2}>t^{1/2}$.  
Then asymptotically
\begin{align*}
  &\prob[x](\tau_{\{m,m'\}}>t)\\
    &\le\exp\left\{
         -t^{\eps^+_1/2}
           \left(
             \frac{K_3}{\log^{2\hat\kappa+M}t}
             -K_2 t^{1/2}
              \exp\left\{
                -\frac{\lambda(\Well(m^-_t))e^{\elevation(\Well(m^+_t))}
                 t^{\eps^+_1/2}}{2}
              \right\}
           \right)
       \right\}\\
    &+\exp\left\{
         -t^{\eps^-_1/2}
           \left(
             \frac{K_3}{\log^{2\hat\kappa+M}t}
             -K_2 t^{1/2}
              \exp\left\{
                -\frac{\lambda(\Well(m^-_t))e^{\elevation(\Well(m^-_t))}
                 t^{\eps^-_1/2}}{2}
              \right\}
           \right)
       \right\}\\
    &\le\exp\left\{
         -K_3 \frac{t^{\eps^+_1/2}}{\log^{2\hat\kappa+M}t}
       \right\}
       +\exp\left\{
         -K_3 \frac{t^{\eps^-_1/2}}{\log^{2\hat\kappa+M}t}
       \right\}\\
    &\le t^{-\eps^+_1}+t^{-\eps^+_1}
\end{align*}
as proposed,
since 
$I^\pm\subset\Well(m^\pm_t)$ implies $\elevation(I^\pm)\le\elevation(\Well(m^\pm_t))$, 
and since \eqref{eq:spectral:gap:elevation} implies $1/t \le \lambda(I^\pm)\exp(\elevation(I^\pm))$ asymptotically, 
which implies $\lambda(I^\pm)\exp(\elevation(I^\pm))t^{\eps^\pm_1} \ge t^{-1+\eps^\pm_1}$ asymptotically,
which implies $K_3/\log^{2\hat\kappa+M}t 
\ge K_2 t^{1/2}\exp(-t^{-1+\eps^\pm_1})
\ge K_2 t^{1/2}\exp(-\lambda(I^\pm)\exp(\elevation(I^\pm))t^{-\eps^\pm_1})$ asymptotically.

\subsection{Proof of Lemma~\ref{lemma:choose:leave}}

This is a classic application of Gambler's Ruin, done before by \citet{Solomon:1975} and \citet{Sinai:1982}.  
We solve it for continuous time setup.  
Our conclusion~\eqref{eq:bound:term2} comes with some straightfoward calculation on the next
\begin{propo}\label{prop:gamblers:ruin}
If $a,z,b\in\ZZ$ are such that $a<z<b$, then
\[
  \prob[z](\tau_a<\tau_b)
    = \frac{\sum_{i=z}^{b-1}e^{V(i)}}
           {\sum_{j=a}^{b-1}e^{V(j)}}.  
\]
\end{propo}

To establish Proposition~\ref{prop:gamblers:ruin} above, we use the Lyapunov function
$f(x)=\sum_{i=a}^{x-1}e^{V(i)-V(a)}$ that renders $f(\xi_t)$ a martingale with respect to $\prob[z]$, as proposed by \citet{Comets:Menshikov:Popov:1998}, and consider the RWRE $\xi^*_t=\xi_{\min(t,\tau_{\{a,b\}})}$ absorbed at the extremes of the interval $[a,b]$, for which trivially $\prob[z](\tau_{a}<\tau_{b})=\prob[z](\tau^*_{a}<\tau^*_{b})$ and $\prob[z](\tau_{b}<\tau_{a})=\prob[z](\tau^*_{b}<\tau^*_{a})$.  
Since $f(\xi_t)$ is a martingale and $\min(t,\tau_{\{a,b\}})$ is a bounded stopping time, we have $\esper[z](f(\xi^*_t))=\esper[z](f(\xi_0))=f(z)$.  
Besides, $f(\xi^*_t)$ is a bounded martingale and, thus, uniformly integrable, so Optional Stopping Theorem render
$f(z)=\lim_{t\to\infty}\esper[z](f(\xi^*_t))=f(a)\prob[z](\tau^*_a<\tau^*_b)+f(b)\prob[z](\tau^*_b<\tau^*_a)$, which implies, as proposed, that
\[
\prob[z](\tau^*_a<\tau^*_b)
  =\frac{f(b)-f(z)}{f(b)-f(a)}
  =\frac{\sum_{i=z}^{b-1}e^{V(i)-V(a)}}{\sum_{j=a}^{b-1}e^{V(j)-V(a)}}
  =\frac{\sum_{i=z}^{b-1}e^{V(i)}}{\sum_{i=a}^{b-1}e^{V(j)}}.
\]

Now, assume that $W(h^-_t)>W(h^+_t)$ (the proof for the case $W(h^-_t)<W(h^+_t)$ is analogous).
We first consider that
$\prob(\tau_{\{m^-_t,m^+_t\}}=\tau_{m^-_t})=\prob(\tau_{m^-_t}<\tau_{m^+_t})$ and then use Proposition~\ref{prop:gamblers:ruin} with $a=m^-_t$, $b=m^+_t$, and $z=0$ to get
\begin{align*}
\prob&(\tau_{m^-_t}<\tau_{m^+_t})
  =\frac{\sum_{i=0}^{m^+_t-1}e^{V(i)}}{\sum_{i=m^-_t}^{m^+_t-1}e^{V(j)}}
  \le\frac{m^+_t\cdot\exp\left(\max_{i=0,\ldots,m^+_t}V(i)\right)}
          {\exp\left(\max_{j=m^-_t,\ldots,m^+_t}V(j)\right)}\\
  &\le\max(m^+_t,m^-_t)\cdot
      \exp\left(\max_{i=0,\ldots,m^+_t}V(i)-\max_{j=m^-_t,\ldots,m^+_t}V(j)\right)\\
  &\le\max(m^+_t,m^-_t)\cdot
      \exp\left(W(h^+_t)-W(h^-_t)+2\max_{j=m^-_t,\ldots,m^+_t}|V(j)-W(j)|\right)\\
  &\le\log^M t\cdot
      \exp\left(-|W(h^+_t)-W(h^-_t)|+2\hat\kappa M\log\log t\right)\\
  &\le t^{-|W(h^+_t)-W(h^-_t)|/\log t}\log^{(2\hat\kappa+1)M} t
\end{align*}
as proposed, since, by hypothesis, $\omega\in\Typical[1][t]\cap\Typical[2][t]$.

\subsection{Proof of Lemma~\ref{lemma:stay:well}}

In our case, Lemma 3.4 from \citet{Comets:Popov:2003} gives $\prob[m](\tau_h<s)\le K_1(s+1)e^{-V(h)+V(m)}$ for every $s\in(0,t]$, which implies
\begin{align*}
\prob[m](\tau_{\{h,h'\}}<t) 
	&\le\prob[m](\tau_h<t)+\prob[m](\tau_{h'}<t) \\
	&\le K_1(t+1)e^{-V(h)+V(m)}+K_1(t+1)e^{-V(h')+V(m)} \\
	&\le K_2\; t\; e^{-\min\{V(h),V(h')\}+V(m)}\\
	&\le K_2\;t\;\exp\left(-\depth(\Well(m))+2\max_{x=h,m,h'}|V(x)-W(x)|\right) \\
	&\le K_2\;t^{-\left(\frac{\depth(\Well(m))}{\log t}-1\right)} \log^{2\hat\kappa M} t
\end{align*}
as proposed, since $\omega\in\Typical[1][t]\cap\Typical[2][t]$.

\subsection{Proof of Lemma~\ref{lemma:oscillating:bottom}}\label{sec:proof:lemma:oscillating:bottom}

We use the reflected version $\hat\xi$ of the RWRE $\xi$ in an interval ($\Well(m)$ in this case) defined in Section~\ref{sec:reflected:rw} above.  

For $s<t$ and $J:=\Well(m)\smallsetminus\Neigh(m)$
\begin{multline*}
\prob[m](\xi_s\not\in\Neigh(m)|\tau_{\{h,h'\}}>s)
	=\prob[m][\hat\omega](\hat\xi_s\not\in\Neigh(m)|\hat\tau_{\{h,h'\}}>s)\\
	=\sum_{x\in J}\prob[m][\hat\omega](\hat\xi_s=x|\hat\tau_{\{h,h'\}}>s)
	=\sum_{x\in J}\prob[m][\hat\omega](\hat\xi_s=x).
\end{multline*}

For $x\in\Well(m)\smallsetminus\Neigh(m)$, the reversibility of $\hat\xi$ and the definition of $\Neigh(m)$ give
\begin{multline*}
\prob[m][\hat\omega](\hat\xi_s=x)
	\le \theta_x/\theta_m 
	\le K_1 e^{-V(x)+V(m)} \\
	\le K_1 e^{-W(x)+W(m)+2\max_{y=x,m}|V(y)-W(y)|}
	\le K_1 t^{-\eps} e^{2\max_{y=x,m}|V(y)-W(y)|}.  
\end{multline*}

By hypothesis, $\omega\in\Typical[1][t]\cap\Typical[2][t]$ gives
\begin{multline*}
\prob[m][\hat\omega](\xi_s\not\in\Neigh(m)|\tau_{\{h,h'\}}>s)
  \le\sum_{x\in J} K_1 t^{-\eps} e^{2\max_{y=x,m}|V(y)-W(y)|}\\
  \le\sum_{x\in J} K_1 t^{-\eps} \log^{2\hat\kappa M} t
  \le K_1 t^{-\eps} \log^{2\hat\kappa M} t |\Well(m)|
  \le K_2 t^{-\eps} \log^{(2\hat\kappa+1)M} t
\end{multline*}
as proposed, since $|\Well(m)|\le|[h^{--}_t,h^{++}_t]|\le 2\log^M t$.


\end{document}

%% file: tstability2.pdftex_t

\begin{picture}(0,0)%
\includegraphics{tstability2.pdf}%
\end{picture}%
\setlength{\unitlength}{3947sp}%
\begingroup\makeatletter\ifx\SetFigFont\undefined%
\gdef\SetFigFont#1#2#3#4#5{%
  \reset@font\fontsize{#1}{#2pt}%
  \fontfamily{#3}\fontseries{#4}\fontshape{#5}%
  \selectfont}%
\fi\endgroup%
\begin{picture}(6818,3490)(112,-2748)
\put(382,-1212){\makebox(0,0)[b]{\smash{{\SetFigFont{9}{10.8}{\familydefault}{\mddefault}{\updefault}{\color[rgb]{0,0,0}$h^{--}_t$}%
}}}}
\put(1389,-938){\makebox(0,0)[b]{\smash{{\SetFigFont{9}{10.8}{\familydefault}{\mddefault}{\updefault}{\color[rgb]{0,0,0}$m^-_t$}%
}}}}
\put(4230,-1266){\makebox(0,0)[b]{\smash{{\SetFigFont{9}{10.8}{\familydefault}{\mddefault}{\updefault}{\color[rgb]{0,0,0}$h^+_t$}%
}}}}
\put(2757,-1266){\makebox(0,0)[b]{\smash{{\SetFigFont{9}{10.8}{\familydefault}{\mddefault}{\updefault}{\color[rgb]{0,0,0}$h^-_t$}%
}}}}
\put(5154,-950){\makebox(0,0)[b]{\smash{{\SetFigFont{9}{10.8}{\familydefault}{\mddefault}{\updefault}{\color[rgb]{0,0,0}$m^+_t=m_t$}%
}}}}
\put(4659,-475){\makebox(0,0)[lb]{\smash{{\SetFigFont{9}{10.8}{\familydefault}{\mddefault}{\updefault}{\color[rgb]{0,0,0}$\eps_2\log t$}%
}}}}
\put(2507,-1718){\makebox(0,0)[lb]{\smash{{\SetFigFont{9}{10.8}{\familydefault}{\mddefault}{\updefault}{\color[rgb]{0,0,0}$(1+\eps^-_3)\log t$}%
}}}}
\put(1940,-700){\makebox(0,0)[rb]{\smash{{\SetFigFont{9}{10.8}{\familydefault}{\mddefault}{\updefault}{\color[rgb]{0,0,0}$(1-\eps^+_1)\log t$}%
}}}}
\put(1500,658){\makebox(0,0)[b]{\smash{{\SetFigFont{9}{10.8}{\familydefault}{\mddefault}{\updefault}{\color[rgb]{0,0,0}$C \log^2 t$}%
}}}}
\put(6165,-2053){\makebox(0,0)[lb]{\smash{{\SetFigFont{9}{10.8}{\familydefault}{\mddefault}{\updefault}{\color[rgb]{0,0,0}$(1+\eps^+_3)\log t$}%
}}}}
\put(6630,-1212){\makebox(0,0)[b]{\smash{{\SetFigFont{9}{10.8}{\familydefault}{\mddefault}{\updefault}{\color[rgb]{0,0,0}$h^{++}_t$}%
}}}}
\end{picture}%

%% file: wellbottom.pdftex_t

\begin{picture}(0,0)%
\includegraphics{wellbottom.pdf}%
\end{picture}%
\setlength{\unitlength}{3947sp}%
\begingroup\makeatletter\ifx\SetFigFont\undefined%
\gdef\SetFigFont#1#2#3#4#5{%
  \reset@font\fontsize{#1}{#2pt}%
  \fontfamily{#3}\fontseries{#4}\fontshape{#5}%
  \selectfont}%
\fi\endgroup%
\begin{picture}(3128,3622)(858,-2771)
\put(2272,-2722){\makebox(0,0)[b]{\smash{{\SetFigFont{11}{13.2}{\familydefault}{\mddefault}{\updefault}{\color[rgb]{0,0,0}$D_{\eps\log t}(m)$}%
}}}}
\put(1116,-2257){\makebox(0,0)[b]{\smash{{\SetFigFont{11}{13.2}{\familydefault}{\mddefault}{\updefault}{\color[rgb]{0,0,0}$h$}%
}}}}
\put(2205,-2257){\makebox(0,0)[b]{\smash{{\SetFigFont{11}{13.2}{\familydefault}{\mddefault}{\updefault}{\color[rgb]{0,0,0}$m$}%
}}}}
\put(3516,-2257){\makebox(0,0)[b]{\smash{{\SetFigFont{11}{13.2}{\familydefault}{\mddefault}{\updefault}{\color[rgb]{0,0,0}$h'$}%
}}}}
\put(2272,529){\makebox(0,0)[b]{\smash{{\SetFigFont{11}{13.2}{\familydefault}{\mddefault}{\updefault}{\color[rgb]{0,0,0}$C_1\log^2 t$}%
}}}}
\put(858,-1484){\makebox(0,0)[rb]{\smash{{\SetFigFont{11}{13.2}{\familydefault}{\mddefault}{\updefault}{\color[rgb]{0,0,0}$\eps\log t$}%
}}}}
\put(3986,-900){\makebox(0,0)[lb]{\smash{{\SetFigFont{11}{13.2}{\familydefault}{\mddefault}{\updefault}{\color[rgb]{0,0,0}$\log t$}%
}}}}
\end{picture}%

%% file: elev.pdftex_t

\begin{picture}(0,0)%
\includegraphics{elev.pdf}%
\end{picture}%
\setlength{\unitlength}{3947sp}%
\begingroup\makeatletter\ifx\SetFigFont\undefined%
\gdef\SetFigFont#1#2#3#4#5{%
  \reset@font\fontsize{#1}{#2pt}%
  \fontfamily{#3}\fontseries{#4}\fontshape{#5}%
  \selectfont}%
\fi\endgroup%
\begin{picture}(5362,3420)(171,-2569)
\put(721,599){\makebox(0,0)[b]{\smash{{\SetFigFont{12}{14.4}{\familydefault}{\mddefault}{\updefault}$f(x)$}}}}
\put(334,-1827){\makebox(0,0)[b]{\smash{{\SetFigFont{12}{14.4}{\familydefault}{\mddefault}{\updefault}$a$}}}}
\put(5327,-1820){\makebox(0,0)[b]{\smash{{\SetFigFont{12}{14.4}{\familydefault}{\mddefault}{\updefault}$b$}}}}
\put(3361,-1827){\makebox(0,0)[b]{\smash{{\SetFigFont{12}{14.4}{\familydefault}{\mddefault}{\updefault}$z$}}}}
\put(4321,-1381){\makebox(0,0)[b]{\smash{{\SetFigFont{12}{14.4}{\familydefault}{\mddefault}{\updefault}$y$}}}}
\put(2521,-1321){\makebox(0,0)[b]{\smash{{\SetFigFont{12}{14.4}{\familydefault}{\mddefault}{\updefault}$x$}}}}
\put(948,-1267){\makebox(0,0)[rb]{\smash{{\SetFigFont{12}{14.4}{\familydefault}{\mddefault}{\updefault}$\elevation(I)$}}}}
\end{picture}%

%% file: newproof.bbl
\begin{thebibliography}{13}
\providecommand{\natexlab}[1]{#1}
\providecommand{\url}[1]{\texttt{#1}}
\expandafter\ifx\csname urlstyle\endcsname\relax
  \providecommand{\doi}[1]{doi: #1}\else
  \providecommand{\doi}{doi: \begingroup \urlstyle{rm}\Url}\fi

\bibitem[Andreoletti(2005)]{Andreoletti:2005}
Pierre Andreoletti.
\newblock Alternative proof for the localization of {S}inai's walk.
\newblock \emph{J. Stat. Phys.}, 118\penalty0 (5--6):\penalty0 883--933, 2005.
\newblock \doi{10.1007/s10955-004-2122-x}.

\bibitem[Andreoletti(2006)]{Andreoletti:2006}
Pierre Andreoletti.
\newblock On the concentration of {S}inai's walk.
\newblock \emph{Stochastic Process. Appl.}, 116\penalty0 (10):\penalty0
  1377--1408, 2006.
\newblock \doi{10.1016/j.spa.2004.12.008}.

\bibitem[Andreoletti(2007)]{Andreoletti:2007}
Pierre Andreoletti.
\newblock Almost sure estimates for the concentration neighborhood of
  {S}inai’s walk.
\newblock \emph{Stochastic Process. Appl.}, 117\penalty0 (10):\penalty0
  1473--1490, 2007.
\newblock \doi{10.1016/j.spa.2007.02.002}.

\bibitem[Berkes et~al.(2014)Berkes, Liu, and Wu]{Berkes:Liu:Wu:2014}
István Berkes, Weidong Liu, and Wei~Biao Wu.
\newblock Komlós–major–tusnády approximation under dependence.
\newblock \emph{The Annals of Probability}, 42\penalty0 (2):\penalty0 794--817,
  03 2014.
\newblock \doi{10.1214/13-AOP850}.
\newblock URL \url{http://dx.doi.org/10.1214/13-AOP850}.

\bibitem[Brox(1986)]{Brox:1986}
Th. Brox.
\newblock A one-dimensional diffusion process in a {W}iener medium.
\newblock \emph{Ann. Probab.}, 14\penalty0 (4):\penalty0 1206--1218, 1986.

\bibitem[Comets and Popov(2003)]{Comets:Popov:2003}
Francis Comets and Serguei Popov.
\newblock Limit law for transition probabilities and moderated deviations for
  {S}inai's random walk in random environment.
\newblock \emph{Probab. Theory Relat. Fields}, 126:\penalty0 571--609, 2003.

\bibitem[Comets et~al.(1998)Comets, Menshikov, and
  Popov]{Comets:Menshikov:Popov:1998}
Francis Comets, Mikhail Menshikov, and Serguei Popov.
\newblock {L}yapunov functions for random walks and strings in random
  environment.
\newblock \emph{Ann. Probab.}, 26\penalty0 (4):\penalty0 1433--1445, 1998.

\bibitem[Koml\'{o}s et~al.(1975)Koml\'{o}s, Major, and
  Tusn\'{a}dy]{Komlos:Major:Tusnady:1975}
J\'{a}nos Koml\'{o}s, P\'{e}ter Major, and G\'{a}bor Tusn\'{a}dy.
\newblock An approximation of partial sums of independent {RV}'s and the sample
  {DF}. {I}.
\newblock \emph{Zeit. Wahrsch. verw. Geb.}, 32:\penalty0 111--131, 1975.

\bibitem[Koml\'{o}s et~al.(1976)Koml\'{o}s, Major, and
  Tusn\'{a}dy]{Komlos:Major:Tusnady:1976}
J\'{a}nos Koml\'{o}s, P\'{e}ter Major, and G\'{a}bor Tusn\'{a}dy.
\newblock An approximation of partial sums of independent {RV}'s and the sample
  {DF}. {II}.
\newblock \emph{Zeit. Wahrsch. verw. Geb.}, 34:\penalty0 33--58, 1976.

\bibitem[Mathieu(1994)]{Mathieu:1994}
Pierre Mathieu.
\newblock Zero {W}hite {N}oise {L}imit through {D}irichlet forms, with
  application to diffusions in random medium.
\newblock \emph{Probab. Theory Relat. Fields}, 99:\penalty0 549--580, 1994.

\bibitem[Sinai(1982)]{Sinai:1982}
Ya.~G. Sinai.
\newblock The limiting behavior of one-dimensional random walk in random
  medium.
\newblock \emph{Theory Probab. Appl.}, 27:\penalty0 256--268, 1982.

\bibitem[Solomon(1975)]{Solomon:1975}
Fred Solomon.
\newblock Random walks in random environments.
\newblock \emph{Ann. Probab.}, 3:\penalty0 1--31, 1975.

\bibitem[Zeitouni(2004)]{Zeitouni:RWRE}
Ofer Zeitouni.
\newblock Random walk in random environment.
\newblock In Jean Picard, editor, \emph{Lectures on Probability Theory and
  Statistics. Ecole d'Et\'{e} de Probabilit\'{e} de Saint-Flour XXXI}, volume
  1837 of \emph{Lecture Notes in Mathematics}, pages 190--312. Springer,
  Berlin, 2004.

\end{thebibliography}
